\documentclass[12pt]{article}
\usepackage{amsmath,amssymb,theorem}
\usepackage{graphicx, enumerate}
\usepackage{subfigure}
\usepackage{psfrag}
\usepackage{placeins}
\usepackage{wrapfig}
\usepackage{booktabs}
\usepackage{chngcntr}
\usepackage{epstopdf}
\counterwithin{figure}{section}
\numberwithin{figure}{section}
\counterwithin{table}{section}
\usepackage{hyperref}
\hypersetup{
  colorlinks,
  citecolor=black,
  filecolor=black,
  linkcolor=black,
  urlcolor=black
}
\newtheorem{thm}{Theorem}[section]
\newtheorem{conj}[thm]{Conjecture}
\newtheorem{cor}[thm]{Corollary}
\newtheorem{lem}[thm]{Lemma}
\newtheorem{prop}[thm]{Proposition}
\theorembodyfont{\rmfamily}
\def\pf{\bigskip\noindent {\bf Proof.}~~}

\def\dfn#1{{\sl #1}}
\def\mytextindent#1{\indent\llap{#1\enspace}\ignorespaces}

\def\myitem{\par\hangindent\parindent\mytextindent}

\newcommand{\less}{\backslash}
\def\proofsquare{  \bigskip\hfill\vrule height3pt width6pt depth2pt}
\def\mc{\mathcal}
\newcounter{counter}
\def\pr#1{(\ref{#1})}

\newcounter{SectionStep}[section]

\begin{document}
\title{Saturation numbers for Ramsey-minimal graphs}
\author{Martin Rolek
~and~ Zi-Xia Song\thanks{Corresponding author. \newline
 E-mail addresses:  msrolek@wm.edu (M. Rolek), Zixia.Song@ucf.edu (Z-X. Song). }\\
Department of Mathematics\\
University of Central Florida\\
 Orlando, FL 32816
}

\maketitle
\begin{abstract}

Given graphs $H_1, \dots, H_t$, a graph $G$ is \dfn{$(H_1, \dots, H_t)$-Ramsey-minimal} if every $t$-coloring of the edges of $G$ contains a monochromatic $H_i$  in color $i$ for some  $i\in\{1, \dots, t\}$, but any proper subgraph of  $G $ does not possess this property. We define  $\mathcal{R}_{\min}(H_1, \dots, H_t)$ to be the family of $(H_1, \dots, H_t)$-Ramsey-minimal graphs.  A graph $G$ is \dfn{$\mathcal{R}_{\min}(H_1, \dots, H_t)$-saturated} if  no element  of $\mathcal{R}_{\min}(H_1, \dots, H_t)$ is a subgraph of $G$, but for any edge $e$ in $\overline{G}$,  some element of $\mathcal{R}_{\min}(H_1, \dots, H_t)$ is  a subgraph of $G + e$. We define $sat(n, \mathcal{R}_{\min}(H_1, \dots, H_t))$ to be  the minimum number of edges over all $\mathcal{R}_{\min}(H_1, \dots, H_t)$-saturated graphs on $n$ vertices. In 1987, Hanson and Toft conjectured   that $sat(n, \mathcal{R}_{\min}(K_{k_1}, \dots, K_{k_t}) )=  (r - 2)(n - r + 2)+\binom{r - 2}{2} $ for $n \ge r$, where $r=r(K_{k_1}, \dots, K_{k_t})$ is the classical Ramsey number for complete graphs.
The first non-trivial  case of Hanson and Toft's conjecture for sufficiently large $n$ was setteled in 2011, and is so far the only settled case.  Motivated by Hanson and Toft's conjecture,   we study the minimum number of edges over all  $\mathcal{R}_{\min}(K_3, \mathcal{T}_k)$-saturated graphs on $n$ vertices, where $\mathcal{T}_k$ is the family of all trees on $k$ vertices.  We  show  that  for $n \ge 18$, $sat(n, \mathcal{R}_{\min}(K_3, \mathcal{T}_4)) =\lfloor {5n}/{2}\rfloor$.  For $k \ge 5$ and $n \ge 2k + (\lceil k/2 \rceil +1) \lceil k/2 \rceil -2$, we obtain an asymptotic bound for $sat(n, \mathcal{R}_{\min}(K_3, \mathcal{T}_k))$ by showing  that  $\left( \frac{3}{2} + \frac{1}{2} \left\lceil \frac{k}{2} \right\rceil \right) n -c\le sat(n, \mathcal{R}_{\min}(K_3, \mathcal{T}_k)) \le \left( \frac{3}{2} + \frac{1}{2} \left\lceil \frac{k}{2} \right\rceil \right) n + C$, where  $c=\left(\frac{1}{2} \left\lceil \frac{k}{2} \right\rceil + \frac{3}{2} \right) k -2$ and $C=  2k^2-6k+\frac32-\left\lceil\frac{k}{2}\right\rceil \left(k- \frac{1}{2}\left\lceil\frac{k}{2}\right\rceil -1\right)$. 
\end{abstract}


\baselineskip=18pt

\section{Introduction}\label{sec:Intro}

All graphs considered in this paper are finite and without loops or multiple edges. For a graph $G$, we will use $V(G)$ to denote the vertex set, $E(G)$ the edge set, $|G|$ the number of vertices, $e(G)$ the number of edges,  $\delta(G)$ the minimum degree, $\Delta(G)$ the maximum degree,  and $\overline{G}$ the complement of $G$.
Given vertex sets $A, B \subseteq V(G)$, we say that $A$ is \dfn{complete to} (resp. \dfn{anti-complete to}) $B$ if for every $a \in A$ and every $b \in B$, $ab \in E(G)$ (resp. $ab \notin E(G)$).
The subgraph of $G$ induced by $A$, denoted $G[A]$, is the graph with vertex set $A$ and edge set $\{xy \in E(G): x, y \in A\}$. We denote by $B \less A$ the set $B - A$, $e_G(A, B)$ the number of edges between $A$ and $B$ in $G$, and $G \less A$ the subgraph of $G$ induced on $V(G) \less A$, respectively.
If $A = \{a\}$, we simply write $B \less a$, $e_G(a, B)$, and $G \less a$, respectively.  For  any  edge  $e\in E(\overline{G})$, we use $G+e$ to denote the graph obtained from $G$ by adding the new edge $e$. 
The {\dfn{join}} $G\vee H$ (resp.~{\dfn{union}} $G\cup H$) of two 
vertex disjoint graphs
$G$ and $H$ is the graph having vertex set $V(G)\cup V(H)$  and edge set $E(G)
\cup E(H)\cup \{xy:   x\in V(G),  y\in V(H)\}$ (resp. $E(G)\cup E(H)$).
Given two isomorphic graphs $G$ and $H$, we may (with a slight but common abuse of notation) write $G = H$.   For an integer $t\ge1$ and a graph $H$, we define $tH$ to be the union of $t$ disjoint copies of $H$.   We use $K_n$, $K_{1,{n-1}}$, $C_n$, $P_n$ and $T_n$ to denote the complete graph,  star,  cycle,   path and a tree on  $n$ vertices, respectively.  \medskip

Given  graphs $G$, ${H}_1, \dots, {H}_t$, we write \dfn{$G \rightarrow ({H}_1, \dots, {H}_t)$} if every $t$-edge-coloring of $G$ contains a monochromatic  ${H}_i$ in color $i$ for some $i\in\{1,2, \dots, t\}$.
The classical \dfn{Ramsey number}  $r({H}_1, \dots, {H}_t)$ is the minimum positive integer $n$ such that $K_n \rightarrow ({H}_1, \dots, {H}_t)$.
A graph $G$ is \dfn{$({H}_1, \dots, {H}_t)$-Ramsey-minimal} if $G \rightarrow ({H}_1, \dots, {H}_t)$, but for any proper subgraph $G'$ of $G$, $G' \not\rightarrow ({H}_1, \dots, {H}_t)$.
We define  $\mathcal{R}_{\min}({H}_1, \dots, {H}_t)$ to be  the family of $({H}_1, \dots, {H}_t)$-Ramsey-minimal graphs. It is straightforward to prove by induction that a graph $G$  satisfies  $G \rightarrow ({H}_1, \dots, {H}_t)$ if and only if there exists a subgraph $G'$ of $G$ such that $G'$ is $({H}_1, \dots, {H}_t)$-Ramsey-minimal.  Ramsey's theorem~\cite{Ramsey1930} implies that $\mathcal{R}_{\min}({H}_1, \dots, {H}_t)\ne\emptyset$ for all integers $t$ and all finite graphs $H_1, \dots, H_t$. As pointed out in a recent paper of Fox, Grinshpun, Liebenau, Person, and Szab\'o~\cite{Fox2016}, ``it is still widely open to classify the graphs in $\mathcal{R}_{\min}(H_1, \dots, H_t)$, or even to prove that these graphs have certain properties".  
Some properties of $\mathcal{R}_{\min}({H}_1, \dots, {H}_t)$  have been studied, such as 
the minimum  degree $s({H}_1, \dots, {H}_t) := \min\{\delta(G) : G \in \mathcal{R}_{\min}({H}_1, \dots, {H}_t)\}$, which was first introduced by Burr, Erd\H os,  and Lov\'asz~\cite{Burr1976}. Recent  results on $s({H}_1, \dots, {H}_t) $ can be found in \cite{Fox2007, Fox2016}. 
For more information on Ramsey-related topics, the readers are referred to  a very recent informative survey due to  Conlon, Fox, and Sudakov~\cite{Conlon2015}.\medskip

In this paper, we study the following problem.  
A graph $G$ is \dfn{$\mathcal{R}_{\min}(H_1, \dots, H_t)$-saturated} if  no element  of $\mathcal{R}_{\min}(H_1, \dots, H_t)$ is a subgraph of $G$, but for any edge $e$ in $\overline{G}$,  some element of $\mathcal{R}_{\min}(H_1, \dots, H_t)$ is  a subgraph of $G + e$. This notion was initiated by Ne\v{s}et\v{r}il~\cite{Nesetril1986} in 1986 when he asked whether there are infinitely many $\mathcal{R}_{\min}(H_1, \dots, H_t)$-saturated graphs. This was answered in the positive by Galluccio, Siminovits, and Simonyi~\cite{Galluccio1992}. 
We define  $sat(n, \mathcal{R}_{\min}(H_1, \dots, H_t))$ to be  the minimum number of edges over all $\mathcal{R}_{\min}(H_1, \dots, H_t)$-saturated graphs on $n$ vertices. This notion was first discussed by  Hanson and Toft~\cite{Hanson1987} in 1987 when  $H_1, \dots, H_t$ are complete graphs. 
They proposed the following conjecture.

\begin{conj}\label{HTC} 
Let $r = r(K_{k_1}, \dots, K_{k_t})$ be the classical Ramsey number for complete graphs. Then
\[ sat(n, \mathcal{R}_{\min}(K_{k_1}, \dots, K_{k_t})) = \displaystyle\left\{
\begin{array}{ll}
\binom{n}{2} \,\, & n < r \\[10pt]

 (r - 2)(n - r + 2) + \binom{r - 2}{2} \,\, & n \ge r
\end{array} \right. \]
\end{conj}
\medskip

Chen, Ferrara, Gould, Magnant, and Schmitt~\cite{Chen2011} proved that  $sat(n, \mathcal{R}_{\min}(K_3, K_3)) = 4n - 10$ for $n\ge56$. This settles the first non-trivial   case of Conjecture~\ref{HTC}  for sufficiently large $n$, and is so far the only settled case.   Ferrara, Kim, and Yeager~\cite{Ferrara2014} proved that  $sat(n, \mathcal{R}_{\min}(m_1K_2, \dots, m_tK_2))=3(m_1+\cdots+m_t-t)$  for $m_1, \dots, m_t\ge1$ and $n>3(m_1+\cdots+m_t-t)$. The problem of finding $sat(n, \mathcal{R}_{\min}(K_3, T_k))$ was also explored in~\cite{Chen2011}.

\begin{prop}\label{prop:satRminKtTm}
Let  $k\ge2$ and  $t\ge2$ be  integers.  Then
\begin{align*}
sat(n, \mathcal{R}_{\min}(K_t, T_k)) \le n(t -& 2)(k - 1) - (t - 2)^2(k - 1)^2 + \binom{(t - 2)(k - 1)}{2} \\
 &+ \left\lfloor \frac{n}{k - 1} \right\rfloor \binom{k - 1}{2} + \binom{r}{2},
\end{align*}
where $r = n ~ ($\emph{mod} $k - 1)$.
\end{prop}
 
It was conjectured in \cite{Chen2011} that the upper bound in Proposition~\ref{prop:satRminKtTm} is  asymptotically correct. Note that there is only one tree on three vertices, namely, $P_3$.  A slightly better result was obtained for $\mathcal{R}_{\min}(K_3, P_3)$-saturated graphs in \cite{Chen2011}. 

\begin{thm}\label{K3P3}
For $n \ge 11$, $sat(n, \mathcal{R}_{\min}(K_3, P_3)) = \left\lfloor \dfrac{5n}{2} \right\rfloor - 5$.
\end{thm}

 Motivated by Conjecture~\ref{HTC},  we study the following problem. Let $\mathcal{T}_k$ be the family of all trees on $k$ vertices.
Instead of fixing a tree on $k$ vertices as in Proposition~\ref{prop:satRminKtTm},  we  will investigate $sat(n, \mathcal{R}_{\min}(K_3, \mathcal{T}_k))$, where a graph $G$ is \dfn{$(K_3, \mathcal{T}_k)$-Ramsey-minimal} if for any $2$-coloring   $c : E(G) \to \{\text{red, blue} \}$, $G$ has either a red $K_3$ or a blue tree $T_k\in \mathcal{T}_k$, and we define $\mathcal{R}_{\min}(K_3, \mathcal{T}_k)$ to be  the family of  $(K_3, \mathcal{T}_k)$-Ramsey-minimal graphs. 
By Theorem~\ref{K3P3}, we see  that $sat(n, \mathcal{R}_{\min}(K_3, \mathcal{T}_3)) = \lfloor {5n}/{2} \rfloor - 5$ for $n \ge 11$.
In this paper, we prove the following two main results. We  first establish  the exact bound for $sat(n, \mathcal{R}_{\min}(K_3, \mathcal{T}_4))$ for $n\ge18$,  and then obtain  an asymptotic bound for    $sat(n, \mathcal{R}_{\min}(K_3, \mathcal{T}_k))$ for all $k \ge 5$ and   $n \ge 2k + (\lceil k/2 \rceil +1) \lceil k/2 \rceil +2$. 

\begin{thm}\label{K3T4}
For $n \ge 18$,  $sat(n, \mathcal{R}_{\min}(K_3, \mathcal{T}_4)) =\left\lfloor \dfrac{5n}{2}\right\rfloor$.

\end{thm}

\begin{thm}\label{K3Tk}
For any integers  $k \ge 5$ and $n \ge 2k + (\lceil k/2 \rceil +1) \lceil k/2 \rceil -2$, there exist constants $c=\left(\frac{1}{2} \left\lceil \frac{k}{2} \right\rceil + \frac{3}{2} \right) k -2$ and $C=  2k^2-6k+\frac32-\left\lceil\frac{k}{2}\right\rceil \left(k- \frac{1}{2} \left\lceil\frac{k}{2}\right\rceil -1\right)$ such that
\[ \left( \frac{3}{2} + \frac{1}{2} \left\lceil \frac{k}{2} \right\rceil \right)n - c \le  sat(n, \mathcal{R}_{\min}(K_3, \mathcal{T}_k)) \le \left( \frac{3}{2} + \frac{1}{2} \left\lceil \frac{k}{2} \right\rceil \right)n + C.
\]
\end{thm}
\medskip

The constants $c$ and $C$ in Theorem~\ref{K3Tk} are both quadratic in $k$.
We believe that the true value of $sat(n, \mathcal{R}_{\min}(K_3, \mathcal{T}_k))$ is closer to the upper bound in Theorem~\ref{K3Tk}.  To establish the desired lower and upper bounds for each of Theorem~\ref{K3T4} and Theorem~\ref{K3Tk}, we need to  introduce more notation and prove a useful lemma (see Lemma~\ref{blue} below). Given a graph $H$, a graph $G$ is \dfn{$H$-free} if $G$ does not contain $H$ as a subgraph.  For a graph $G$,  let   $c : E(G) \to \{\text{red, blue} \}$  be a  $2$-edge-coloring of $G$ and let  $E_r$ and $E_b$ be the  color classes of the coloring $c$.  We use $G_{r}$ and $G_{b}$ to denote the spanning subgraphs of $G$  with edge sets  $E_r$ and $E_b$, respectively. We define  $c$ to be a  \dfn{bad $2$-coloring} of $G$ if $G$ has neither  a red $K_3$ nor a blue $T_k\in  \mathcal{T}_k$, that is, if  $G_r$ is $K_3$-free and $G_b$ is $T_k$-free for any $T_k\in\mathcal{T}_k$.  For any $v\in V(G)$, we use $d_r(v)$ and $N_r(v) $ to denote the degree and neighborhood of $v$ in $G_r$, respectively. Similarly, we define $d_b(v)$ and $N_b(v)$ to be the degree and neighborhood of $v$ in $G_b$, respectively.  \bigskip

\noindent {\bf Remark.}  One can see that if  $G$ is $\mathcal{R}_{\min}(K_3, \mathcal{T}_k)$-saturated,   then $G$ admits at least one bad $2$-coloring but,  for any edge $e\in E(\overline{G})$, 
$G+e$ admits  no bad $2$-coloring.   \medskip

 We will utilize the following Lemma~\ref{blue}(a) to  force  a unique bad $2$-coloring of certain graphs  in order to establish an   upper bound for $sat(n, \mathcal{R}_{\min}(K_3, \mathcal{T}_k))$.  Lemma~\ref{blue}(b) and  Lemma~\ref{blue}(c) will be applied to establish a lower bound for $sat(n, \mathcal{R}_{\min}(K_3, \mathcal{T}_k))$.

\begin{lem}\label{blue}
For any integer $k\ge3$, let   $c : E(G) \to \{\text{red, blue} \}$  be a bad $2$-coloring of a  graph  $G$ on $n\ge k+2$ vertices. 

\myitem{(a)}  If $e \in E(G)$ belongs to at least $2k - 3$ triangles in $G$, then $e\in E_b$.
\myitem{(b)} If $G$ is  $\mathcal{R}_{\min}(K_3, \mathcal{T}_k)$-saturated and $D_1, \dots, D_p$ are the components of $G_b$ with     $|D_i|< {k}/{2}$ for all $i\in\{1,\dots, p\}$, then $p\le2$.  Moreover, if $p=2$, then  $V(D_1)$ is complete to $V(D_2)$ in $G_r$. 
\myitem{(c)}  If $G$ is  $\mathcal{R}_{\min}(K_3, \mathcal{T}_k)$-saturated, and among all bad $2$-colorings of $G$, $c$ is chosen so that  $|E_r|$ is maximum, then  $\Delta(G_r) \le n-3$ and $G_r$ is 2-connected.
\end{lem}

\pf  To prove (a), suppose that there exists an edge $e = uv \in E_r$ such that $e$ belongs to at least $2k-3$ triangles in $G$.
Since $G_r$ is  $K_3$-free,  we see that either $d_b(u)\ge k-1$ or $d_b(v)\ge k-1$. In either case,   $G_b$ contains $K_{1, {k-1}}$  as a subgraph, a contradiction. \medskip

To prove (b), let  $D_1, \dots, D_p$  be given as in (b). 
We next show that $p\le2$.  Since $G$ is $\mathcal{R}_{\min}(K_3, \mathcal{T}_k)$-saturated, we see that,  for any edge $e$ in $\overline{G}$, $G+e$ admits no bad $2$-coloring.  We claim that,  for any  $i,j\in\{1,\dots, p\}$ with $i\ne j$,   $V(D_i)$ is complete to $V(D_j)$  in $G_r$.  
Suppose 
that there exist vertices   $u \in V(D_i)$ and  $v \in V(D_j)$ such  that  $uv \notin E_r$.  Then $uv\notin E(G)$ and so we obtain a bad $2$-coloring of  $G+uv$ from $c$ by  coloring the edge $uv$ blue, a contradiction. Thus  $V(D_i)$ is complete to $V(D_j)$ in $G_r$ for any  $i,j\in\{1,\dots, p\}$ with $i\ne j$.
Since $G_r$ is $K_3$-saturated, it follows that $p\le2$.\medskip

It remains to prove (c).  By the choice of $c$, $G_{r}$ is $K_3$-free but $G_r+e$ contains a $K_3$ for any $e\in E(\overline{G_r})$,  and $G_b$ is ${T}_k$-free for any $T_k\in \mathcal{T}_k$. 
Note that $G_b$ is disconnected and every component of $G_b$ contains at most $k - 1$ vertices. 
Since $G$ is $\mathcal{R}_{\min}(K_3, \mathcal{T}_k)$-saturated, we see that, for any edge $e$ in $\overline{G}$,  $G+e$ admits no bad $2$-coloring.  
Suppose that   $\Delta(G_r) \ge n-2$. Let  $x \in V(G)$ with $d_{r}(x) = \Delta(G_r)$ and  let $v$ be the  unique non-neighbor of $x$ in $G_r$ if $d_r(x)=n-2$.    
Since $G_r$ is $K_3$-free, we see that $N_r(x)$ is an independent set in $G_r$.  By the choice of $c$,  $v$ must be  complete to $N_r(x)$ in $G_r$. Since $n\ge k+2$, we have $|N_r(x)|\ge k$.  Let $u\in N_r(x)$ and let $H$ be the component of $G_b$ containing $u$. Then $|H|\le k-1$ and $V(H)\subset N_r(x)$. Let $w\in N_r(x)\less V(H)$. Clearly,  $uw\notin E(G)$.  We obtain a bad $2$-coloring of  $G+uw$  from  $c$ by coloring the edge $uw$ red,  and  then recoloring all  edges incident with $u$ in $G_r$  blue and all  edges incident with $u$ in $G_b$ red,    a contradiction.  This proves that $\Delta(G_r) \le n-3$. \medskip

Finally, we show that $G_r$ is $2$-connected. 
Suppose that $G_r$ is not $2$-connected. Since $G_{r}$ is $K_3$-free but $G_r+e$ contains a $K_3$ for any $e\in E(\overline{G_r})$, we see that  $G_r$ is connected and  must have  a cut vertex, say $u$. Since  $\Delta(G_r) \le n-3$, $u$ has a non-neighbor, say  $v$,  in $G_r$. Let $G_1$ and $G_2$ be two components of $G_r\less u$ with  $v \in V(G_2)$.  
Let $w \in V(G_1)$. By the choice of $c$, $wv\notin E_b$, otherwise we obtain a bad $2$-coloring of $G$ from $c$ by recoloring the blue edge $wv$ red. Thus 
 $wv\notin E(G)$ and then we obtain a bad $2$-coloring of  $G+wv$ from $c$ by coloring the edge  $wv$ red,    a contradiction. Therefore $G_r$ is $2$-connected. \medskip

This completes the proof of Lemma~\ref{blue}. \proofsquare

The remainder of this paper is organized as follows. In Section~\ref{sec:K3Sat},  we discuss  $K_3$-saturated graphs with a specified minimum degree and prove a structural result which we shall use in the proof of Theorem~\ref{K3T4}.    We then prove Theorem~\ref{K3T4} in  Section~\ref{sec:satRminK3T4} and Theorem~\ref{K3Tk} in   Section~\ref{sec:K3Tk}.

\section{$K_3$-saturated graphs}\label{sec:K3Sat}

In this section we list  known results and establish  new ones on $K_3$-saturated graphs that we shall need to prove our main results. \medskip

Given a graph $H$,  
a graph $G$ is \dfn{$H$-saturated} if $G$ is  $H$-free but,  for any edge $e \in E(\overline{G})$,   $G + e$ contains a copy of $H$  as a subgraph. We define  $sat(n, H)$ to be  the minimum number of edges over all $H$-saturated graphs on $n$ vertices. This notion
 was introduced by Erd\H os, Hajnal,  and Moon~\cite{Erdos1964} in 1964.  Results on $H$-saturated graphs can be found in  surveys by  either Faudree, Faudree, and Schmitt~\cite{Faudree2011} or Pikhurko~\cite{Pikhurko}.  In this section we are interested in the case when $H=K_t$.
Erd\H os, Hajnal,  and Moon~\cite{Erdos1964} showed that if $G$ is a $K_t$-saturated graph on $n$ vertices, then 
$e(G)\ge(t-2)n- \binom{t - 1}{2}$. Moreover, they showed that the graph $K_{t- 2} \vee  \overline{K}_{n - t + 2}$ is the unique $K_t$-saturated graph with $n$ vertices and $(t-2)n- \binom{t - 1}{2} $ edges. Notice that this extremal graph has minimum degree $t-2$.  One may  ask:  what is the minimum number of edges in a $K_t$-saturated graph with specified minimum degree? This was first studied by Duffus and Hanson~\cite{Duffus1986}  in 1986. They proved the following two results.

\begin{thm}\label{delta=2}
If $G$ is a $K_3$-saturated graph on  $n \ge 5$ vertices with  $\delta(G) = 2$, then $e(G) \ge 2n - 5$ edges. Moreover,  if $e(G) = 2n - 5$,  then $G$ can be  obtained from $C_5$ by repeatedly duplicating vertices of degree $2$. 
\end{thm}

\begin{thm}\label{delta=3}
If $G$ is a $K_3$-saturated graph on  $n \ge 10$ vertices with $\delta(G) = 3$, then $e(G) \ge 3n - 15$. Moreover, if  $e(G) = 3n - 15$,   then $G$ contains the Petersen graph as a subgraph.
\end{thm}

Alon,  Erd\H{o}s,   Holzman,  and Krivelevich~\cite{Alon1996} showed that  any $K_4$-saturated graph  on  $n \ge 11$ vertices with minimum degree $ 4$ has at least  $4n-19$ edges.  This has recently been generalized by Bosse, the second author, and Zhang~\cite{Bosse2017+} by showing that  any $K_t$-saturated graph on  $n \ge t+7$ vertices with minimum degree $ t\ge3$ has at least  $tn-{{t+1}\choose 2}-9$ edges.    Moreover, they showed that the graphs $K_{t- 3} \vee H$ are the only  $K_t$-saturated graphs with $n$ vertices and $tn- \binom{t +1}{2} -9$ edges, where $H$ is a $K_3$-saturated graph on $n-t+3\ge10$ vertices with $\delta(H)=3$.    Theorem~\ref{Kp} below is a  result of  Day~\cite{Day2017} on $K_t$-saturated graphs with prescribed minimum degree. It confirms  a conjecture of Bollob\'as~\cite{Bollobas}  when $t=3$. It is worth noting that the constant $c$ given in Theorem~\ref{Kp} does not have a dependency on $t$. This is a consequence of the fact that  every  $K_t$-saturated graph has minimum degree at least $t-2$.

\begin{thm}\label{Kp}
For any integers $p \ge 1$ and $t \ge 3$, there exists a constant $c = c(p)$ such that if $G$ is a $K_t$-saturated graph on $n$ vertices with $\delta(G) \ge p$, then $e(G) \ge pn - c$.
\end{thm}

\begin{figure}[tbhp]
\centering
\includegraphics[width=200px]{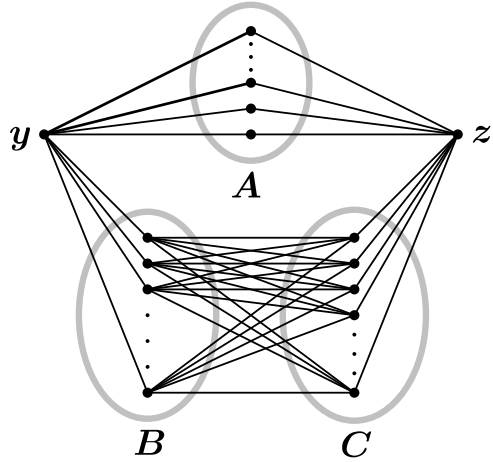}
\caption{The graph $J$ 
}
\label{J}
\end{figure}

 For our proof of Theorem~\ref{K3T4}, we will need a structural result on $K_3$-saturated graphs with minimum degree  at most $ 2$. The graph $J$ depicted  in  Figure~\ref{J}  is a $K_3$-saturated graph with minimum degree $2$, where $A\ne\emptyset$ and  either  $B=C=\emptyset$ or  $B\ne\emptyset$ and  $C \ne \emptyset$;   $A$, $B$ and  $C$ are  independent sets in $J$ and pairwise disjoint;    $A$ is anti-complete to $B\cup C$ and  $B$ is complete to $C$;  $N_J(y)=A\cup B$ and $N_J(z)=A\cup C$;  and $|A|+|B|+|C|=|J|-2$.  It is straightforward to check that 
$e(J) = 2(|J| - 2) + |B||C| - |B| - |C|\ge 2|J|-5$. Moreover, $e(J)=2|J|-5$ when $|B|=1$ or $|C|=1$. That is, $e(J)=2|J|-5$ when  $J$ is  obtained from $C_5$ by   repeatedly duplicating vertices of degree $2$. 
Lemma~\ref{structural} below yields a new proof of Theorem~\ref{delta=2}, and has  been generalized for all $K_t$-saturated graphs with minimum degree at most $ t - 1$ in~\cite{Bosse2017+}.

\begin{lem}\label{structural}
Let $G$ be a $K_3$-saturated graph with  $n$ vertices and  $\delta(G)=\delta$.

\myitem{(a)} If $\delta = 1$, then $G=K_{1, n - 1}$.

\myitem{(b)} If $\delta = 2$, then $G =J$,   where the graph $J$ is  depicted in Figure~\ref{J}.  Moreover, $J=K_{2, n-2}$ when $B=C=\emptyset$. 

\myitem {(c)}	 If $\delta \ge 3$, then $2e(G) \ge \max\{(\delta+1)n-\delta^2-1,\, (\delta+2)n - \delta(\delta + t) -2 \}$, where   $t := \min \{ d(v) : v  \text{ is adjacent to }
\text{a vertex of degree }\delta \text{ in } G\}$.
\end{lem}

\pf Let  $x \in V(G)$ be a vertex with $d(x) = \delta$.
Since $G$ is $K_3$-saturated, we see that  $G$ is connected and $K_3$-free. 
First assume that $d(x) = 1$.  Let $y$ be the neighbor of $x$. If there exists a vertex $z\in V(G)$ such that $yz\notin E(G)$, then $G+xz$ is $K_3$-free, contrary to the fact that $G$ is $K_3$-saturated. 
Thus  $y$ is complete to $V(G) \setminus \{y\}$.
Clearly, $N(y)$ is an independent set because $G$ is $K_3$-free. Thus  $G = K_{1, n-1}$.
This proves (a).\medskip 

Next assume that $d(x) = 2$.   Let $N(x) = \{y, z\}$.
Then $yz \notin E(G)$ because $G$ is $K_3$-free. We next show that $N(y) \cup N(z) = V(G)\less\{y,z\}$.  Suppose there exists a vertex $w\in V(G)$ such that $wy, wz\notin E(G)$. Then $G+xw$ is $K_3$-free, contrary to the fact that $G$ is $K_3$-saturated. 
  Hence   $N(y) \cup N(z) = V(G)\less\{y,z\}$.
Let  $A := N(y) \cap N(z)$, $B := N(y) \setminus N(z)$, and $C := N(z) \setminus N(y)$.  Then  $|A| + |B| + |C| = n - 2$, and $A, B, C$ are pairwise disjoint. 
Clearly,  $x \in A$, and either $B=C=\emptyset$ or $B\ne\emptyset$ and  $C \ne \emptyset$ because $\delta(G)= 2$.   
Since $G$ is   $K_3$-free, we see that $A, B, C$  are independent sets in $G$, and $A$ is anti-complete to $B \cup C$.
We next show that  $B$ must be complete to $C$ when  $B\ne\emptyset$ and  $C \ne \emptyset$. Suppose there exist vertices $b\in B$ and $c\in C$ such that $bc\notin E(G)$. Then $G+bc$ is $K_3$-free, a contradiction. 
Thus $G=J$, where $J$ is depicted in Figure~\ref{J}.  \medskip

It remains to prove (c).   Let $\delta\ge3$ and let $t$ be given as in (c). Then   $d(x)\ge3$.  We first show that $2e(G) \ge (\delta+1)n-\delta^2-1$. Since $G$ is $K_3$-saturated, every vertex in $V(G) \backslash N[x]$ has at least one neighbor in $N(x)$, yielding $ \sum_{v \in N(x)} d(v) \ge |V(G) \backslash N[x]|+d(x)=n-1$. Therefore 
\begin{align*}
2e(G) & =  d(x) + \sum_{v \in N(x)} d(v) + \sum_{v \in V(G) \backslash N[x]} d(v)\\
& \ge  \delta+n-1+\delta(n-\delta-1)\\
& \ge (\delta+1)n-\delta^2-1.
\end{align*}

We next show that $2e(G) \ge  (\delta+2)n - \delta(\delta + t) -2$.  We may assume that there exists a vertex $y \in N(x)$ with $d(y) = t$.   Notice that $x$ and $y$ have no common neighbor. Let $M: =V(G) \backslash ( N(x) \cup N(y) )$. Then $|M|=n-\delta-t$. Since $G$ is $K_3$-saturated, each vertex in  $M$ has at least one neighbor in $N(x) \backslash y$ and at least one neighbor in $N(y) \backslash x$. Thus  $\sum_{v \in N(x) \backslash y} d(v) \ge n-t-1$, and $\sum_{v \in N(y) \backslash x} d(v) \ge n-\delta -1$. Then

\begin{align*}
2e(G) & = d(x)+d(y) + \sum_{v \in N(x) \backslash y} d(v) + \sum_{v \in N(y) \backslash x} d(v) + \sum_{v \in M} d(v) \\
& \ge \delta +t+ (n-t-1) + (n-\delta-1) + \delta(n-\delta-t) \\
& =(\delta+2)n - \delta(\delta + t) -2. 
\end{align*}

This completes the proof of Lemma~\ref{structural}.
\proofsquare

\begin{cor}\label{2n}
Let $G$ be a $K_3$-saturated graph on $n\ge5$ vertices with $\delta(G)=2$.
If $e(G) = 2n - k$ for some $k \in \{0, 1, 2, 3, 4, 5\}$, then  $G = J$ with $|B||C| - |B| - |C| = 4 - k$, where $A, B, C,$ and $J$ are as depicted in Figure~\ref{J} and the values of $|B|$ and $|C|$ are summarized in Table~\ref{table}.
\end{cor}

\pf 
Since $\delta(G) = 2$, by Lemma~\ref{structural}(b), $G =J$ with $e(G) = 2(n - 2) + |B||C| - |B| - |C|$ and either $B = C = \emptyset$ or $B, C \ne \emptyset$, where $A, B, C,$ and $J$ are as depicted in Figure~\ref{J}.
We see that  $|B||C| - |B| - |C| = 4 - k$ because  $e(G) = 2n - k$, where $k \in \{0, 1, 2, 3, 4, 5\}$. 
Solving the resulting equation in each case of $k$ yields explicit constructions of $J$, which are summarized in Table~\ref{table}.
\proofsquare

\begin{table}[htbp]
\centering
\begin{tabular}{ *5l @{} *9l @{}*9l @{} }    \toprule
$k$ & $e(J)$ 	& \emph{values of $|B|$ and $|C|$ with $|B|\le|C|$}   \\\midrule
$5$ & $2n-5$   	&  $|B|=1$ and $|C| \ge 1$ \\ 
$4$ & $2n-4$		& $|B|=|C|=2$ or $|B|=|C|=0$\\ 
 $3$ & $2n-3$		& $|B| =2$ and  $|C| = 3$ \\
$2$ & $2n-2$		& $|B| = 2$ and $|C| = 4$ \\
$1$ & $2n-1$		& $|B|=2$ and $|C|=5$ or $|B|=|C|=3$ \\
$0$ & $2n$		& $|B|=2$ and $|C|=6$ \\\bottomrule
 \hline
\end{tabular}
\caption{Construction of the graph  $J$ determined by $k$}
\label{table}
\end{table}

\section{Proof of Theorem~\ref{K3T4}}\label{sec:satRminK3T4}

\begin{figure}[htbp]
\centering
\includegraphics[scale=.6]{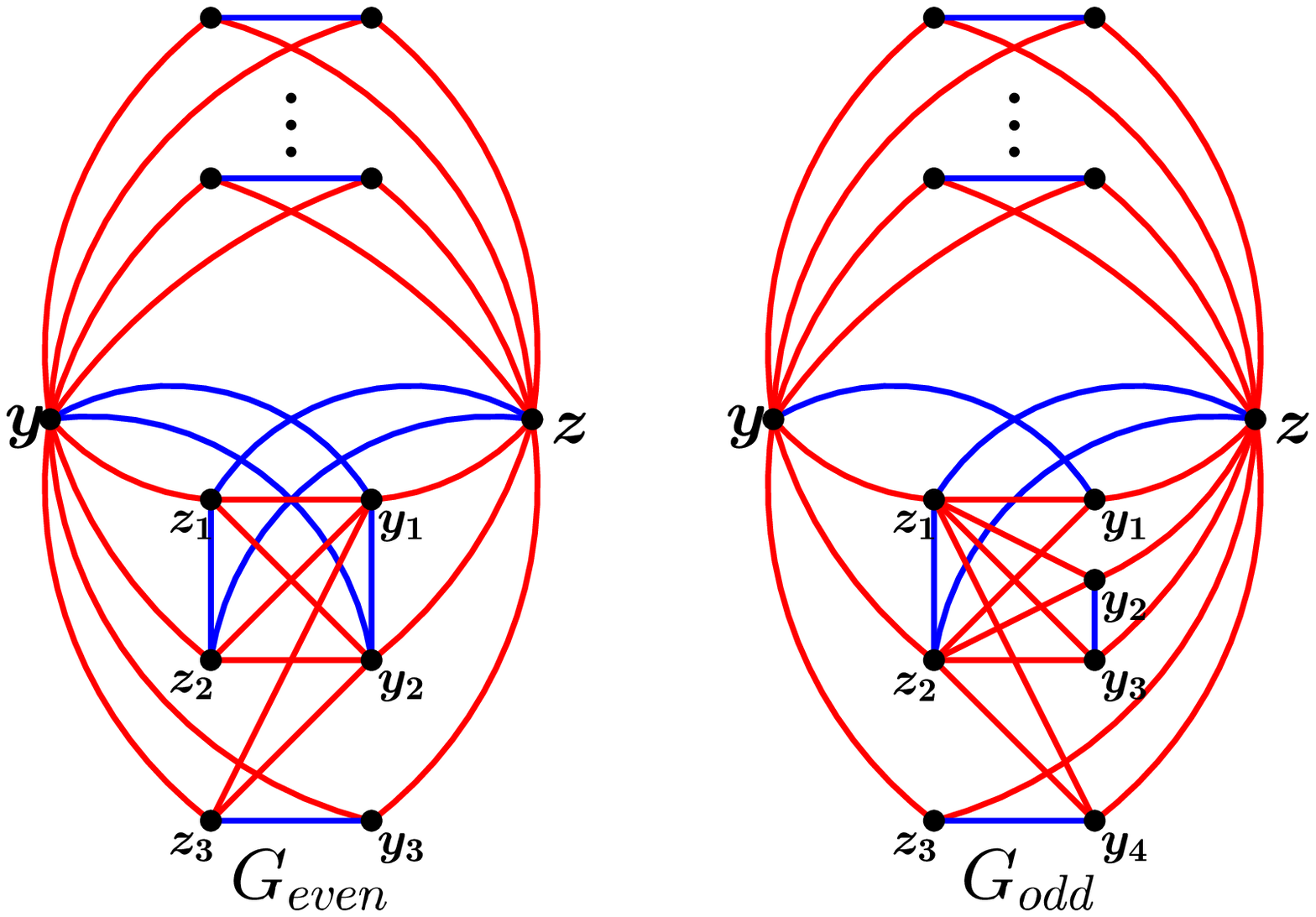}
\caption{Two  $\mc{R}_{\min} (K_3, \mathcal{T}_4)$-saturated graphs with a unique bad  $2$-coloring, where  dashed lines indicate blue and solid lines indicate red.}
\label{EO}
\end{figure}

We are now ready to prove Theorem~\ref{K3T4}.   We first  establish the desired upper bound for $sat(n, \mc{R}_{\min} (K_3, \mathcal{T}_4))$ by    constructing an $\mathcal{R}_{\min}(K_3, \mathcal{T}_4)$-saturated graph with the desired number of edges.  Let $n\ge8 $ be an integer and let $H=(\lfloor n/2\rfloor -4)K_2$.      When $n\ge8$ is even,  let   $G_{even}$  be the graph  obtained from  $H$  by adding eight new vertices $y, z, y_1, y_2, y_3,  z_1, z_2, z_3$,   and then  joining:  $y$ to all vertices in $V(H)\cup\{y_1, y_2, y_3, z_1, z_2, z_3\}$;  $z$ to all vertices in $V(H)\cup\{y_1, y_2, y_3, z_1, z_2\}$;  $y_1$ to all vertices in $\{y_2, z_1,  z_2, z_3\}$;  $y_2$ to all vertices in $\{z_1,  z_2, z_3\}$, $z_1$ to $z_2$;   and  $z_3$ to $y_3$.  When $n$ is odd, let   $G_{odd}$  be the graph  obtained from  $H$  by adding nine new vertices $y, z, y_1, y_2, y_3, y_4, z_1, z_2, z_3$,   and then  joining:  $y$ to all vertices in $V(H)\cup\{y_1,  z_1, z_2, z_3\}$;   $z$ to all vertices in $V(H)\cup\{y_1, y_2, y_3, y_4, z_1, z_2, z_3\}$;  $z_1$ to all vertices in $\{ y_1, y_2, y_3, y_4, z_2\}$;  $z_2$ to all vertices in $\{y_1, y_2, y_3, y_4\}$, $y_2$ to $y_3$;  and $y_4$ to $z_3$.    The graphs $G_{odd}$ and $G_{even}$  are depicted in Figure~\ref{EO}.  It can be easily checked that $e(G_{odd})= (5n-1)/{2}$ and $e(G_{even})= {5n}/{2}$.   We next show   that   $G_{odd}$ and $G_{even}$  are $\mathcal{R}_{\min}(K_3, \mathcal{T}_4)$-saturated.\medskip

One can easily check that the coloring $c : E(G) \to \{\text{red, blue} \}$ for each of $G_{odd}$ and $G_{even}$  given in Figure~\ref{EO} is a bad $2$-coloring.  We next show that $c$ is the unique bad $2$-coloring  for each of  $G_{odd}$ and $G_{even}$.  To find a bad $2$-coloring for  $G_{odd}$,  
by Lemma~\ref{blue}(a), the edges  $zz_1, zz_2, z_1z_2$  must be colored blue and so all the other edges incident with $z, z_1, z_2$ must be red.  Then $yy_1, y_2y_3, y_4z_3$ and all edges in $E(H)$ must be  blue and all the other edges incident with $y$ must be red.  This proves that $G_{odd}$ has a unique bad $2$-coloring, as depicted in Figure~\ref{EO}.  To find a bad $2$-coloring for $G_{even}$,  by Lemma~\ref{blue}(a),   $y_1y_2$ must be colored blue.  We next show that  $z_1z_2$ must be colored blue. Suppose that $z_1z_2$ is colored red. To avoid a red $K_3$, we may assume that $yz_1$ is colored blue. Then all edges $z_1y_1, z_1y_2, yy_1, yy_2$ must be red, and so $z_2y_1, z_2y_2$  must be blue, which then forces $y_1z$ to be red and $z_1z$ to be blue.  Now the edges $z_3y$ and  $z_3y_1$ must be colored red, which yields a red $K_3$ with vertices $y, z_3, y_1$.   This proves that  $z_1z_2$ must be colored blue. Similar to the argument for $G_{odd}$,  one can see that  the coloring of $G_{even}$, depicted in Figure~\ref{EO}, is  the unique bad $2$-coloring of $G_{even}$. 
It is straightforward to see that  both $G_{odd}$ and $G_{even}$ are $\mathcal{R}_{\min}(K_3, \mathcal{T}_4)$-saturated, and so  $sat(n, \mathcal{R}_{\min}(K_3, \mathcal{T}_4)) \le \lfloor {5n}/{2}\rfloor$.  We next show that $sat(n, \mathcal{R}_{\min}(K_3, \mathcal{T}_4)) \ge \lfloor {5n}/{2}\rfloor$.  \medskip

 Let  $G$ be  an $\mc{R}_{\min} (K_3, \mathcal{T}_4)$-saturated graph on $n\ge18$ vertices.  Then, for any edge   $e\in E(\overline{G})$, $G+e$   has no bad $2$-coloring.  Suppose that $e(G)< 5n/2$ if $n$ is even and $e(G)<(5n-1)/2$ if $n$ is odd. 
   Among  all bad $2$-colorings of $G$, let  $c : E(G) \to \{\text{red, blue} \}$ be a bad $2$-coloring of  $G$  with $|E_r|$  maximum.    By the choice of $c$,  $G_r$ is  $K_3$-saturated.  Note that $G_b$ is disconnected and every component of $G$ is isomorphic to $K_1$, $K_2$,  $P_3$ or $K_3$. By Lemma~\ref{blue}(c), we have 
 \bigskip

\noindent \refstepcounter{counter}\label{e:maxdeg} (\arabic{counter})\,\, $\Delta(G_r) \le n-3$ and  $G_r$ is 2-connected.\medskip

We next show that \medskip

\noindent \refstepcounter{counter}\label{e:Gr=J} (\arabic{counter})\,\, 
$\delta(G_r)=2$ and so $G_r=J$ with $A\ne\emptyset$, $B\ne\emptyset$,  and $C\ne\emptyset$, where $J$, $A, B, C$ are  depicted in Figure~\ref{J}.  

\pf By \pr{e:maxdeg}, $\delta(G_r)\ge2$. Suppose that $\delta(G_r)\ge3$.  We next show that $e(G_r)\ge \lceil (5n-17)/2\rceil$. This is trivially true if $\delta(G)\ge5$. So we may assume that $3\le \delta(G_r)\le4$.  By Theorem~\ref{delta=3} applied to $G_r$ when $\delta(G_r)=3$  and Lemma~\ref{structural}(c) applied to $G_r$ when $\delta(G_r)=4$, we see that $e(G_r)\ge \lceil (5n-17)/2\rceil$ because  $n\ge18$.  By Lemma~\ref{blue}(b), $e(G_b)\ge \lceil (n-2)/2\rceil$. Thus $e(G)=e(G_r)+e(G_b) \ge \lfloor {5n}/{2}\rfloor$, a contradiction. Hence  $\delta(G_r)=2$. By Lemma~\ref{structural}(b), $G_r=J$, where $J$, $A\ne\emptyset, B, C$ are  depicted in Figure~\ref{J}.  By \pr{e:maxdeg}, $B\ne\emptyset$ and $C\ne\emptyset$.  \proofsquare

 For the remainder of the proof, let $J$, $A, B, C$,   and $y, z$ be  given as  in Figure~\ref{J}, where $A\ne\emptyset$, $B\ne\emptyset$,  and $C\ne\emptyset$.   By \pr{e:Gr=J}, $G_r=J$.  We next show that \bigskip
 
 \noindent \refstepcounter{counter}\label{e:bc=2} (\arabic{counter})\,\, 
$|B|\ge 2$ and  $|C| \ge 2$.

\pf  Suppose that $|B|=1$ or $|C| = 1$, say the latter.  Let $u$ be the vertex in $C$.   If $yz, yu\in E_b$, then $d_b(u)=1$ because $G_b$ is $T_4$-free. Now for any  $w\in A$, we obtain a bad $2$-coloring of  $G+uw$  from  $c$ by coloring the edge $uw$ red,  and  then recoloring the edge $zu$ blue. Thus  either $yz\notin E(G)$ or $yu\notin E(G)$. We may assume that $yz\notin E(G)$.  Then $yu\in E_b$, otherwise, we obtain a bad $2$-coloring of  $G+yu$  from  $c$ by coloring the edge $yu$ blue, and then recoloring the edge $zu$ blue, and all the edges incident with $z$ and $u$ in $G_b$ red.
Notice that  $d_b(u) = 1$, for otherwise  let $w \in A$ be the other neighbor of $u$ in $G_b$ and $v \in B$. Then $d_b(w)=1$ and so we obtain a bad $2$-coloring of  $G+wv$  from  $c$ by coloring the edge $wv$ red,  and then recoloring the edge $yw$ blue.  We next claim that $B=N_b(z)$.\medskip

Suppose that $B \ne N_b(z)$.
Let $w \in B \setminus N_b(z)$, and let $K$ be the component of $G_b$ containing $w$.
If $V(K)\subseteq B$, then for any  $v\in A$, we obtain a bad $2$-coloring of  
$G+wv$  from  $c$ by coloring the edge $wv$ red,  and  then  recoloring the edges $yw, uw$ blue and all edges incident with $w$ in $G_b$ red, a contradiction. Thus $V(K)\cap A\ne\emptyset$.  Let $v\in V(K)\cap A$. We claim that $V(K)=\{w, v\}$. Suppose that $|K|=3$. Let $v'$ be the third vertex of $K$. Then $K$ is isomorphic to $K_3$. If $v'\in A$, then we obtain a bad $2$-coloring of $G$ from $c$ by recoloring the edge $yw$ blue, and then recoloring the edges $wv, wv'$ red, contrary to the choice of $c$. Thus $v'\in B$, which again  yields a bad $2$-coloring of $G$ from $c$ by recoloring the edge $yv$ blue,  and then recoloring the edges $vw, vv'$ red, contrary to the choice of $c$. Thus $V(K)=\{w, v\}$, as claimed. 
For any $v^* \in (A\cup B)\less (\{w,v\}\cup N_b(z))$, we obtain a bad $2$-coloring of $G+wv^*$ from $c$ by coloring the edge $wv^*$ red,  and then recoloring the edge $wv$ red, and the edges $yw, uw$ blue.
Thus $B=N_b(z)$, as  claimed.\medskip

Since $B=N_b(z)$, we have $|B|\le2$.  Then $yu\in E_b$, otherwise by a similar argument for showing $B=N_b(z)$, we have $|A|=|N_b(u)|\le2$ and so $n\le 7$, a contradiction.   Let $v\in B$.  If $B = \{v\}$, then by a similar argument  for  showing   $d_b(u)=1$,  we have  $d_b(v) = 1$.
But then  we obtain a bad $2$-coloring of  $G+yz$  from  $c$ by coloring the edge $yz$ blue, and then   recoloring  the edge $yu$ red,  and the edge $yv$ blue. 
Thus  $|B|=2$. Let   $v'$ be the other vertex in $B$. Then  $vv' \in E_b$, otherwise we obtain a bad $2$-coloring of  $G+vv'$  from  $c$ by coloring the edge $vv'$ blue. 
But now we obtain a bad $2$-coloring of  $G+yz$  from  $c$ by coloring the edge $yz$ blue, and then  recoloring the  edges $yu, vv', zv'$ red, and  edges $yv, uv'$ blue, a contradiction.   
\proofsquare

 By Lemma~\ref{blue}(b), $G_b$ has at most two isolated vertices. Thus  $e(G_b)\ge (n-2)/2$. Since $e(G)< 5n/2$, we see that $e(G_r)\le 2n$.  By \pr{e:bc=2}, $|B|\ge2$ and $|C|\ge2$.  By Corollary~\ref{2n},  $e(G_r)\ge 2n-4$ and $|B|+|C|\le8$. Thus $|A|\ge n-10\ge8$. We next show that \bigskip

\noindent \refstepcounter{counter}\label{e:P3} (\arabic{counter})\,\, 
If  $P_3$ is a component of $G_b\less \{y,z\}$ with vertices $x_1, x_2, x_3$ in order, then $x_2\in A$  and $|\{x_1, x_3\}\cap B|=|\{x_1, x_3\}\cap C|=1$.

\pf  Clearly, $\{x_1, x_2, x_3\}\not\subseteq A\cup B$ or $\{x_1, x_2, x_3\}\not\subseteq A\cup C$, otherwise $x_1x_3\notin E(G)$ and we obtain  a bad $2$-coloring of $G+x_1x_3$ from $c$ by coloring the edge  $x_1x_3$ blue.  Since $y, z\notin \{x_1, x_2, x_3\}$, we see that $x_2\in A$. Then $|\{x_1, x_3\}\cap B|=|\{x_1, x_3\}\cap C|=1$.
\proofsquare

\noindent \refstepcounter{counter}\label{e:yz} (\arabic{counter})\,\, 
 $yz\notin E(G)$.  \medskip

\pf Suppose that $yz\in E(G)$. Then $yz\in E_b$.  Since $G_b$ does not contain a  ${T}_4$, we see that either $d_b(y)=1$ or $d_b(z)=1$. We may assume that $d_b(z)=1$.   We claim that $d_b(y)=1$ as well.  Suppose that $d_b(y)=2$. Let  $w\in C$ be the other neighbor of $y$ in $G_b$. Then $d_b(w)=1$. Let $v\in A$.  We obtain a bad $2$-coloring of  $G+wv$  from  $c$ by coloring the edge $wv$ red, and   recoloring the edge $zw$ blue. Thus  $d_b(y)=d_b(z)=1$.  Since $e(G_r)\le 2n$ and $|A|\ge n-10\ge8$, by Corollary~\ref{2n} and \pr{e:P3}, $G_b$ contains a component, say $K$, such that $V(K)\cap A\ne\emptyset$ and $V(K)\subset A\cup B$ or $V(K)\subset A\cup C$.  Let $u\in V(K)\cap A$ and $w\in A\less V(K)$.   We obtain a bad $2$-coloring of  $G+uw$  from  $c$ by coloring the edge $uw$ red, and then  recoloring the  edges $yu, zu$ blue,  and all the edges incident with $u$ in $G_b$  red,   a contradiction.  \proofsquare

\noindent \refstepcounter{counter}\label{e:Gb} (\arabic{counter})\,\, 
$G_b$ has  no  isolated vertex.  

\pf Suppose for a contradiction that $G_b$ has an  isolated vertex, say $u$.   Then $d(u)=d_r(u)$. By \pr{e:maxdeg}, $d(u)\le n-3$.  For any  $w\in V(G) \backslash N[u]$,   adding a blue edge $uw$ to $G$ must  yield a blue ${T}_4$, because  $G$ is $\mc{R}_{\min} (K_3, \mathcal{T}_4)$-saturated.  Hence, \medskip

\noindent ($*$) \,  every  vertex of   $V(G) \backslash N[u]$ belongs  to a $P_3$ or $K_3$ in $G_b$.  \medskip

We next claim that every vertex of $A\less u$ belongs to a $P_3$ or $K_3$ in $G_b$. By ($*$), this is obvious if $u\in A\cup B\cup C$. So we may assume that $u\in \{y,z\}$. By symmetry, we may further assume that $u=z$. By \pr{e:yz}, $yz\notin E(G)$.   Suppose that there exists a vertex $v\in A$ such that $v$ belongs to a component, say $K$, with $|K|\le2$. Then $V(K)\subseteq A\cup B$ or $V(K)\subseteq A\cup C$. Let  $w\notin V(K)$ be a vertex  in $C$. This is possible  because $|C|\ge 2$ by \pr{e:bc=2}. 
We then  obtain a bad $2$-coloring of  $G+vw$  from  $c$ by coloring the edge $vw$ red, and   recoloring   the edge $vu$ blue, a contradiction. Thus every vertex of $A\less u$ belongs to a $P_3$ or $K_3$ in $G_b$, as claimed. \medskip

Since $|B|+|C|\le8$ and $|A|\ge n-10\ge8$, by \pr{e:P3} and Corollary~\ref{2n}, we see that $G_b[A]$ has at least two components isomorphic to   $K_3$. By Lemma~\ref{blue}(b), $G_b$ has at most two isolated vertices and so $e(G_b)\ge 6+(n-8)/2$. Since $e(G)<5n/2$, we have  $e(G_r)\le 2n-3$.  By \pr{e:bc=2}, $|B|\ge2$ and $|C|\ge2$.  By Corollary~\ref{2n},  $2n-4\le e(G_r)\le 2n-3$ and $\max\{|B|,|C|\}\le3$.   Thus $|A|\ge n-8\ge10$.  
By \pr{e:P3} and Corollary~\ref{2n} again,  $G_b[A]$ has at least three components isomorphic to   $K_3$. Thus $e(G_b)\ge 9+\lceil(n-11)/2\rceil$ and so $e(G)\ge (2n-4)+9+\lceil(n-11)/2\rceil\ge\lfloor 5n/2\rfloor$, a contradiction. \proofsquare

\noindent \refstepcounter{counter}\label{e:1nbr} (\arabic{counter})\,\, 
  $d_b(y)=d_b(z)=2$. 

\pf  Suppose that $d_b(y)\le1$ or $d_b(z)\le1$.  By \pr{e:Gb}, $d_b(y), d_b(z)\ge1$.  We may assume that $d_b(y)=1$.   By  \pr{e:yz}, $yz\notin E(G)$.  Let $y_1\in C$ be the unique neighbor of $y$ in $G_b$, and let  $z_1\in B$ be a neighbor of $z$ in $G_b$. We claim that $d_b(y_1)=1$. Suppose that  $d_b(y_1)=2$. Let  $y_1^*\in A\cup C$ be the other neighbor of $y_1$ in $G_b$.  Then $y_1^*\in A$, otherwise,    we obtain a bad $2$-coloring of  $G+yy_1^*$  from  $c$ by coloring the edge   $yy_1^*$  blue.  Let $w\in B$.  Then we obtain a bad $2$-coloring of  $G+y_1^*w$  from  $c$ by coloring the edge $y_1^*w$  red and  recoloring the edge $y_1^*y$ blue.  
Thus  $d_b(y_1)=1$, as claimed.\medskip

 By \pr{e:bc=2}, $ |B|\ge2$ and  $|C|\ge 2$. We next claim that $N_b(z)=B$.  Suppose that there exists a vertex $u\in B$ such that $uz\notin E(G_b)$. Then $uz_1\notin E_b$, otherwise,   we obtain a bad $2$-coloring of  $G+uz$  from  $c$ by coloring the edge   $uz$  blue.  This implies that  $B\less N_b(z)$ is anti-complete to $N_b(z)$ in $G_b$.  Let $K$ be the component of $G_b$ containing $u$. By \pr{e:Gb},  $|K|\ge2$.  Since $G_b$ is $T_4$-free, we see that $N_b[z]$ is anti-complete to $V(K)$ in $G_b$.  
Suppose first that $V(K) \subseteq B$. 
If $K$ is isomorphic to $K_3$ or $|N_b(z)|=2$, then $|B| \ge 4$ and $G_b$ contains at least one $K_3$ ($K$ or $G[N_b[z]]$).
By Corollary~\ref{2n}, $e(G_r) \ge 2n - 2$.  By \pr{e:Gb}, $e(G_b) \ge 3 + \lceil (n - 3)/2 \rceil$.
Hence $e(G) =e(G_r)+e(G_b)\ge (2n - 2) + 3 + \lceil(n - 3)/2\rceil \ge \lfloor 5n/2 \rfloor$, a contradiction.
Thus $K$ is isomorphic to $K_2$ and $d_b(z)=1$. 
Using a similar argument to show that $d_b(y_1)=1$, we have $d_b(z_1) = 1$. Let $V(K)=\{u, u'\}$. If $B=\{u,u', z_1\}$, then   we obtain a bad 2-coloring of $G + yz$ from $c$  by coloring the edge $yz$ blue, and then recoloring the edges $y_1u, y_1u', yz_1$ blue,  and  the edge $yy_1$ red. Thus $|B|\ge4$. By Corollary~\ref{2n}, $|C|=2$. Let $C=\{w, y_1\}$.  Let $v\in A$ be such that $v$ and $w$ are not in the same component of $G_b$.  This is possible because $|A|\ge 8$. Then 
 we obtain a bad 2-coloring of $G + vw$ from $c$ by coloring the edge $vw$ red, and then recoloring the edges $z_1w, zw$ blue, and all the edges incident with $w$ in $G_b$ red. 
This proves that $V(K) \not\subseteq B$ and so $V(K)\cap A\ne\emptyset$.  Let $v\in V(K)\cap A$. We next show  that $V(K)=\{u, v\}$. Suppose that $|K|=3$. Let $v'$ be the third vertex of $K$. Then $K$ is isomorphic to $K_3$. If $v'\in A$, then we obtain a bad $2$-coloring of $G$ from $c$ by recoloring the edge $uy$ blue,  and then recoloring the edges $uv, uv'$ red, contrary to the choice of $c$. If  $v'\in B$, then we obtain  a bad $2$-coloring of $G$ from $c$ by recoloring the edge $vy$ blue,  and then recoloring the edges $vu, vv'$ red, contrary to the choice of $c$. Thus $v'\in C$. Now for any $w\in A\less v$, we obtain  a bad $2$-coloring of $G+uw$ from $c$ by coloring the edge $uw$ red,  and then recoloring the edges $uy,  uy_1$ blue, and  $uv$ red. 
Hence $V(K)=\{u, v\}$. For any $v' \in A\less v$, we obtain a bad $2$-coloring of $G+uv'$ from $c$ by coloring the edge $uv'$ red,  and then recoloring the edges $uy$ blue and  $uv$ red.  
 Thus    $N_b(z)=B$, as claimed.  \medskip
 
Since $N_b(z)=B$  and $d_b(z)\le2\le |B|$, we see that  $|B|=2$. Let $B=\{z_1, z_2\}$. Then $z_1z_2\in E(G_b)$, otherwise, we obtain a bad $2$-coloring of  $G+z_1z_2$  from  $c$ by coloring the edge   $z_1z_2$  blue.    Let $C=\{y_1,  \dots, y_t\}$, where $t=|C|$.  Then $y_1y_j\notin E(G_b)$ for all $j\in\{2, \dots, t\}$ because $d_b(y_1)=1$.  If $t\ge4$, then by Corollary~\ref{2n}, $e(G_r)\ge 2n-2$. By \pr{e:Gb},   $e(G_b)\ge 3+\lceil(n-3)/2\rceil$. Thus  $e(G)\ge (2n-2)+3+\lceil(n-3)/2\rceil\ge\lfloor 5n/2\rfloor$, a contradiction. Thus $2\le t\le 3$.  Let $v\in A$ be such that $vy_j\notin E(G)$ for all $j\in\{1, 2, \dots, t\}$. This is possible because $|A|\ge 8$ and $t\le3$. We obtain a bad $2$-coloring of  $G+y_2v$  from  $c$ by coloring the edge   $y_2v$  red, and then  when $t=2$,   recoloring  the edges $yz_1, z_1y_1, z_2y_2,y_2z$ blue,  the edges $z_1z, z_1z_2$,  and all the edges incident with $y_2$ in $G_b$  red;   when $t=3$,  recoloring  the edges $y_1z_1, y_1z_2, zy_2,zy_3$ blue, the edges $yy_1, zz_1, zz_2$,  and all the edges between $A$ and  $\{y_2, y_3\}$ in $G_b$  red.     \proofsquare

By \pr{e:1nbr}, $d_b(y)=d_b(z)=2$.  By \pr{e:yz}, $yz\notin E(G)$.  Let $N_b(y)=\{y_1, y_2\}\subseteq C$ and $N_b(z)=\{z_1, z_2\}\subseteq B$.  Then $y_1y_2, z_1z_2\in E_b$, otherwise,  we obtain a bad $2$-coloring of  $G+e$  from  $c$ by coloring the edge   $e$  blue, where $e\in \{y_1y_2, z_1z_2\}$.  By \pr{e:Gb},  $e(G_b)\ge 6+\lceil(n-6)/2\rceil$. Since $e(G)< \lfloor 5n/2\rfloor$, by Corollary~\ref{2n}, we see that $n$ is even and $|B|=|C|=2$.  Let $v\in A$.  We obtain a bad $2$-coloring of  $G+vz_1$  from  $c$ by coloring the edge   $vz_1$  red, and then recoloring the edges  $yz_1, z_2y_1,   z_2y_2$ blue,  and  edges $yy_1, yy_2, zz_2,  z_1z_2 $ red,  a contradiction. \medskip

This completes the proof of Theorem~\ref{K3T4}. \hfill\vrule height3pt width6pt depth2pt\\

\section{Proof of Theorem~\ref{K3Tk}}\label{sec:K3Tk}

Finally, we prove Theorem~\ref{K3Tk}.  We will construct an $\mathcal{R}_{\min}(K_3, \mathcal{T}_k)$-saturated graph on $n \ge 2k + (\lceil k/2 \rceil +1)\lceil k/2 \rceil-2$ vertices  which yields the desired   upper bound in Theorem~\ref{K3Tk}.  \medskip

\begin{figure}[htbp]
\centering
\includegraphics[width=400px]{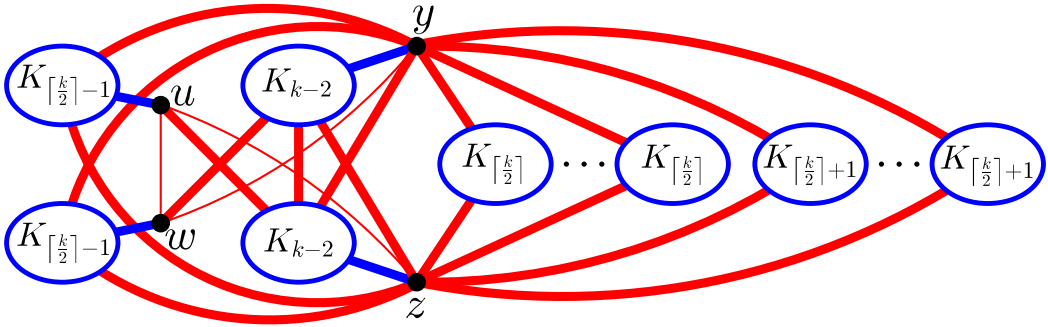}
\caption{An $\mathcal{R}_{\min}(K_3, \mathcal{T}_k)$-saturated graph with a unique bad $2$-coloring, where  dashed lines indicate blue and solid lines indicate red.}
\label{SatK3Tk}
\end{figure}

For positive  integers $k, n$ with  $k\ge 5$ and  $n \ge 2k + (\lceil k/2 \rceil +1)\lceil k/2 \rceil-2$, let $t$ be the remainder of  $n-2k-2\lceil k/2\rceil+2$ when divided by $\lceil k/2 \rceil$,  and  let $H= 2 K_{\lceil k/2\rceil-1}\cup 2K_{k-2} \cup s K_{\lceil k/2\rceil}\cup tK_{\lceil k/2\rceil+1}$, where $s\ge0$ is an integer satisfying   $s\lceil k/2\rceil+ t(\lceil k/2\rceil+1)=n-2k-2\lceil k/2\rceil+2$.  Let $H_1$, $H_2$ be the  two disjoint copies of $K_{k-2}$,   and let $H_3, H_4$ be  the two disjoint   copies of   $K_{\lceil k/2\rceil-1}$ in $H$, respectively. 
Finally,   let $G$ be the  graph obtained from $H$ by adding four new vertices $y, z, u, w$,  and  then joining:  every vertex in $H_1$ to all vertices in $H_2$;  $y$ to all vertices in $V(H)\cup \{w\}$;  $z$ to  all vertices in $V(H)\cup \{u\}$;  $u$ to all vertices in $\{w\}\cup V(H_2)\cup V(H_3)$; and  $w$ to all vertices in $ V(H_1)\cup V(H_4)$,  as depicted in Figure~\ref{SatK3Tk}.

Clearly,  the coloring $c : E(G) \to \{\text{red, blue} \}$  given in Figure~\ref{SatK3Tk} is a bad $2$-coloring of $G$. We next show that  $c$ is the unique bad $2$-coloring of $G$.   By Lemma~\ref{blue}(a), each edge $e\in E(H_1)\cup E(H_2)$ must be colored blue because $e$ belongs to $2k-3$ triangles in $G$. Then all edges between $V(H_1)$ and $V(H_2)$ in $G$ must be colored red and  the edge $yv$ must be colored red for some $v\in V(H_1)\cup V(H_2)$,  because $G_b$ is $T_k$-free.
Additionally, $y$ can only be joined by a blue edge to a vertex in either $V(H_1)$ or $V(H_2)$ but not both. It follows that $y$ is complete to one of $V(H_1)$ or $V(H_2)$ in $G_r$. We next show that $y$ is complete to $V(H_2)$ in $G_r$. 
Suppose that $y$ is complete to $V(H_1)$ in $G_r$.
Then $y$ is complete to $V(H_2)$ in $G_b$ since $G_r$ is $K_3$-free, and so $yw \in E_r$ since $G_b$ is $T_k$-free.
This implies that  $z$ must be   complete to $V(H_1)$ in $G_b$.
But now $w$ must be complete to $V(H_1)$ in $G_r$,  which yields a red $K_3$ on $y, w, v$ for any $v\in V(H_1)$, a contradiction. 
Hence $y$ is complete to $V(H_2)$ in $G_r$.
Then $y$ must be  complete to $V(H_1)$ in $G_b$.
Since $G_b$ is $T_k$-free, $y$ is complete to $\{w\}\cup(V(H)\less V(H_1))$ in $G_r$, and $z$ is complete to $V(H_1)$ in $G_r$. Since $G_r$ is $K_3$-free, we see that all edges in each component of $H$ must be colored blue, and then $z$ must be complete to $V(H_2)$ in $G_b$ and $w$ must be  complete to $V(H_4)$ in $G_b$.  By symmetry of  $y$ and $z$, it follows that $z$ is complete to $\{u\}\cup(V(H)\less V(H_2))$ in $G_r$, and $u$ is complete to $V(H_3)$ in $G_b$.
 This proves that $c$ is the unique bad $2$-coloring of $G$.
It is straightforward to see that  $G$ is $\mathcal{R}_{\min}(K_3, \mathcal{T}_k)$-saturated. Using the facts  that $s\lceil k/2\rceil+ t(\lceil k/2\rceil+1)=n-2k-2\lceil k/2\rceil+2$ and $t\le \lceil k/2\rceil -1$, we see that   
\begin{align*}
e(G) & = 2(n-2)+{2k-4\choose 2} +(2(k-2)+1)+(s+2) {\lceil k/2\rceil\choose 2}+ t{\lceil k/2\rceil+1\choose 2} \\
& =(2n+2k^2-7k+3)+(s+2)\lceil k/2\rceil \frac{\lceil k/2\rceil-1}2+t(\lceil k/2\rceil+1)\frac{(\lceil k/2\rceil-1)+1}2\\
&=(2n+2k^2-7k+3)+ \dfrac{\lceil k/2\rceil-1}2 \left((s+2)\lceil k/2\rceil+t(\lceil k/2\rceil+1)\right) + \dfrac{t}{2} \left(\left\lceil\frac{k}{2}\right\rceil+1\right)\\
&=(2n+2k^2-7k+3)+ \dfrac{\lceil k/2\rceil-1}2 \left((s\lceil k/2\rceil+t(\lceil k/2\rceil+1))+2\lceil k/2\rceil\right) + \dfrac{t}{2} \left(\left\lceil\frac{k}{2}\right\rceil+1\right)\\
&\le (2n+2k^2-7k+3)+ \dfrac{\lceil k/2\rceil-1}2 \left(n-2k-2\lceil k/2\rceil+2+2\lceil k/2\rceil\right) + \dfrac{t}{2} \left(\left\lceil\frac{k}{2}\right\rceil+1\right)\\
 &\le \left( \frac{3}{2} + \frac{1}{2} \left\lceil\frac{k}{2}\right\rceil \right) n + 2k^2-6k+2-(k-1)\left\lceil\frac{k}{2}\right\rceil+ \dfrac{\lceil k/2\rceil-1}{2} \left(\left\lceil\frac{k}{2}\right\rceil+1\right)\\
  &\le \left( \frac{3}{2} + \frac{1}{2} \left\lceil\frac{k}{2}\right\rceil \right) n + 2k^2-6k+\frac32-\left\lceil\frac{k}{2}\right\rceil \left(k- \dfrac{1}{2} \left\lceil\frac{k}{2}\right\rceil -1\right)\\
 &= \left( \frac{3}{2} + \frac{1}{2} \left\lceil\frac{k}{2}\right\rceil \right) n + C, 
\end{align*}

\noindent  where $C=2k^2-6k+\frac32-\left\lceil\frac{k}{2}\right\rceil \left(
k- \frac{1}{2} \left\lceil\frac{k}{2}\right\rceil -1\right)$. Therefore  $sat(n, \mathcal{R}_{\min}(K_3, \mathcal{T}_k)) \le e(G)\le  \left( \frac{3}{2} + \frac{1}{2} \left\lceil\frac{k}{2}\right\rceil \right) n +C$.  \medskip
 
 Let $c=\left(\frac{1}{2} \left\lceil \frac{k}{2} \right\rceil + \frac{3}{2} \right) k -2$. We next show that   $sat(n, \mathcal{R}_{\min}(K_3, \mathcal{T}_k)) \ge \left( \frac{3}{2} + \frac{1}{2} \left\lceil\frac{k}{2}\right\rceil \right) n -c$.   Let $G$ be an $\mathcal{R}_{\min}(K_3, \mathcal{T}_k)$-saturated graph on $n \ge 2k + (\lceil k/2 \rceil +1) \lceil k/2 \rceil -2$ vertices. Then $G+e$   has no bad $2$-coloring  for any edge   $e\in E(\overline{G})$.  
   Among  all bad $2$-colorings of $G$, let  $c : E(G) \to \{\text{red, blue} \}$ be a bad $2$-coloring of  $G$  with $|E_r|$  maximum.     
By the choice of $c$, $G_{r}$ is $K_3$-saturated and $G_b$ is ${T}_k$-free for any $T_k\in \mathcal{T}_k$.
Note that $G_b$ is disconnected and every component  of $G_b$ contains at most $k - 1$ vertices. By Lemma~\ref{blue}(c), we have 
 \bigskip
\setcounter{counter}{0}

\noindent \refstepcounter{counter}\label{e:maxdeg} (\arabic{counter})\,\, $\Delta(G_r) \le n-3$ and  $G_r$ is 2-connected.\medskip

Let $D_1, D_2, \dots, D_p$ be the components of $G_b$. Since $n \ge 2k + (\lceil k/2 \rceil +1) \lceil k/2 \rceil -2$, we have $p\ge 3$. We next show that \medskip

\noindent \refstepcounter{counter}\label{e:clique} (\arabic{counter})  \, $G[V(D_i)] =K_{|D_i|}$ for all $i\in\{1,2, \dots, p\}$. 

\pf Suppose that there exists a component   of $G_b$, say $D_1$,  such that $G[V(D_1)] \ne K_{|D_1|}$. Let  $u, v \in V(D_1)$  be such that $uv \notin E(G)$.
We obtain a bad $2$-coloring of $G + uv$ from $c$ by coloring the edge $uv$ blue, a contradiction. \proofsquare\medskip

\noindent \refstepcounter{counter}\label{e:edgecount} (\arabic{counter})  
\, $ \displaystyle \sum_{i = 1}^p e(G[V(D_i)]) \ge  \left(\frac{1}{2} \left\lceil \frac{k}{2} \right\rceil - \frac{1}{2} \right) n -\left(\frac{1}{2} \left\lceil \frac{k}{2} \right\rceil - \frac{1}{2} \right) k$

\pf  By \pr{e:clique},  $G[V(D_i)] =K_{|D_i|}$ for all $i\in\{1,2, \dots, p\}$.  By Lemma~\ref{blue}(b),  at most two components  $D_i$ have fewer than $k/2$ vertices. Let $t$ be the remainder of  $n - k$ when divided by $\lceil k/2 \rceil$, and let $s \ge 0$ be an integer such that $n - k = s \lceil k/2 \rceil + t (\lceil k/2 \rceil + 1)$. 
It is straightforward to see that $\displaystyle \sum_{i = 1}^p e(G[V(D_i)])$  is minimized when:  two of the components, say $D_1, D_2$, are such that $|D_1|, |D_2|< k/2$;   $t$ of the components, say $D_3, \dots, D_{t+2}$, are such that  $|D_3 |=\cdots=|D_{t+2}| =\lceil k/2 \rceil + 1$;  and  $s$ of the components, say $D_{t+3}, \dots, D_{t+s+2}$,  are such that  $|D_{t+3}|=\cdots=|D_{t+s+2}|=  \lceil k/2 \rceil$.
Using the facts  that $s\lceil k/2\rceil+ t(\lceil k/2\rceil+1)=n-2k-2\lceil k/2\rceil+2$ and $t\le \lceil k/2\rceil -1$,  it  follows that 
\begin{align*}
\sum_{i = 1}^p e(G[V(D_i)])  &> s {\lceil k/2\rceil\choose 2} + t {\lceil k/2\rceil + 1\choose 2} \\
 &= s \lceil k/2\rceil \frac{\lceil k/2\rceil-1}2+t(\lceil k/2\rceil+1)\frac{(\lceil k/2\rceil-1)+1}2\\
 &= \left(\frac{1}{2} \left\lceil \frac{k}{2} \right\rceil - \frac{1}{2} \right)\left(s\left\lceil \frac{k}2\right\rceil+t\left(\left\lceil \frac{k}2\right\rceil+1\right)\right) + \dfrac{t}{2} \left(\left\lceil\frac{k}{2}\right\rceil+1\right)\\
 &\ge \left(\frac{1}{2} \left\lceil \frac{k}{2} \right\rceil - \frac{1}{2} \right)(n-k) \\
 &\ge \left(\frac{1}{2} \left\lceil \frac{k}{2} \right\rceil - \frac{1}{2} \right) n - \left(\frac{1}{2} \left\lceil \frac{k}{2} \right\rceil - \frac{1}{2} \right) k.
\end{align*}
\proofsquare

Assume that $G_b[V(D_i)] =K_{|D_i|}$ for all $i\in\{1,2, \dots, p\}$. By \pr{e:edgecount},  $ |E_b| \ge  \left(\frac{1}{2} \left\lceil \frac{k}{2} \right\rceil - \frac{1}{2} \right) n -\left(\frac{1}{2} \left\lceil \frac{k}{2} \right\rceil - \frac{1}{2} \right) k$.
By Lemma~\ref{blue}(b) and Theorem~\ref{delta=3}, 
$ |E_r| \ge 2n - 5$. 
Therefore $e(G) = |E_r| + |E_b| \ge \left( \frac{3}{2} + \frac{1}{2} \left\lceil \frac{k}{2} \right\rceil \right) n  -\left(\frac{1}{2} \left\lceil \frac{k}{2} \right\rceil - \frac{1}{2} \right) k -5\ge \left( \frac{3}{2} + \frac{1}{2} \left\lceil \frac{k}{2} \right\rceil \right) n -c$, where $c=\left(\frac{1}{2} \left\lceil \frac{k}{2} \right\rceil + \frac{3}{2} \right) k -2$, as desired.  So we may assume that $G_b[V(D_i)] \ne K_{|D_i|}$ for some  $i\in\{1,2, \dots, p\}$, say $i=1$. 
Let  $u_1, u_2 \in V(D_1)$  be such that $u_1u_2\notin E_b$.   By \pr{e:clique}, $u_1u_2\in E_r$.   Since $G_r$ is $K_3$-saturated, we have $N_r(u_1)\cap N_r(u_2)=\emptyset$.  We next show that \bigskip

\noindent \refstepcounter{counter}\label{e:rednbr} (\arabic{counter}) \,\,
for any  $j\in\{2,\dots, p\}$ and any  $w \in V(D_j)$, if $wu_i\notin E_r$ for some $i\in\{1,2\}$, then $N_r(w)\cap N_r(u_i)\less (V(D_1)\cup V(D_j))\ne\emptyset$. 

\pf We may assume that $wu_1\notin E_r$. Since $G_r$ is $K_3$-saturated, we see that $N_r(w)\cap N_r(u_1)\ne\emptyset$. Note that $wu_1\notin E(G)$. If $N_r(w)\cap N_r(u_1)\less (V(D_1)\cup V(D_j))=\emptyset$, then we obtain a bad $2$-coloring of $G+wu_1$ from $c$ by coloring $wu_1$ red,  and then recoloring all red edges incident with $u_1$ in $D_1$ blue and all red edges incident with $w$ in  $D_j$ blue, a contradiction. \proofsquare

\noindent \refstepcounter{counter}\label{e:2rednbr} (\arabic{counter}) \,\, 
For any  $j\in\{2,\dots, p\}$ and any  $w \in V(D_j)$,  $|N_r(w)\less  V(D_j)|\ge2 $. 

\pf   This is obvious when $wu_1, wu_2\in E_r$.  So we may assume that $wu_1\notin E_r$.     Since $N_r(u_1)\cap N_r(u_2)=\emptyset$, it follows from  \pr{e:rednbr} that either 
$|N_r(w)\less (V(D_1)\cup V(D_j))|\ge2$ when $wu_2\notin E(G)$ or   $|N_r(w)\less  V(D_j)|=|N_r(w)\less (V(D_1)\cup V(D_j))|+|N_r(w)\cap V(D_1)|\ge1+1=2 $ when  $wu_2\in E(G)$.   In both cases,  $|N_r(w)\less  V(D_j)|\ge2 $, as desired.  \proofsquare

For each vertex $w \in V(G) \less V(D_1)$, since $G_r$ is $K_3$-saturated, we see that either  $wu_1\notin E_r$ or  $wu_2\notin E_r$.  Let $P :=\{ w \in V(G) \setminus V(D_1): \, wu_1, wu_2\notin E_r\}$, $Q :=\{ w \in V(G) \setminus V(D_1): \, wu_1\notin E_r, wu_2\in E_r  \}$, and $R :=\{ w \in V(G) \setminus V(D_1): \, wu_1\in E_r, wu_2\notin E_r\}$.  Further, let $Q_1$ denote the set of vertices  $w\in Q$ such that $N_r(w)\cap V(D_1)=\{u_2\}$,  and  let $R_1$  denote the set of vertices  $w\in R$ such that $N_r(w)\cap V(D_1)=\{u_1\} $. Let $Q_2:=Q\less Q_1$ and $R_2:=R\less R_1$.  By definition,  $P, Q_1, Q_2, R_1, R_2$ are pairwise disjoint and $|P|+|Q|+|R|=n-|V(D_1)|\ge n-k+1$.   Let $H$ be obtained from $G\less V(D_1)$ by deleting all edges in $G[V(D_i)]$ for all $i\in\{2,3,\dots, p\}$. Then $E(H)\subset E_r$ and for each edge $e$ in $H$, $e$ is not in $G[V(D_i)]$ for any  $i\in\{2,3,\dots, p\}$. For any $w\in Q_1\cup R_1$, by \pr{e:rednbr}, $N_{H}(w)\less P\ne\emptyset$. We next show that \bigskip

\noindent \refstepcounter{counter}\label{e:1rednbr} (\arabic{counter}) \,\, 
for any     $w \in Q_1$,  if $w$ is adjacent to exactly one vertex, say $v$, in   $H\less P$, then $v\in R_2 $. 

\pf We may assume that $w\in V(D_2)$.  Since $w\in Q_1$, we have $N_r(w)\cap V(D_1)=\{u_2\}$.  By \pr{e:rednbr}, $vu_1\in E_r$,  and we may further assume that $v\in V(D_3)$.  Then $vu_2\notin E_r$ because $G_r$ is $K_3$-free. Since $D_1$ is a component of $G_b$, there must exist a vertex, say $u\in V(D_1)$, such that $uu_2\in E_b$. Then $wu\notin E_r$  (and so $wu\notin E(G)$) because $N_r(w)\cap V(D_1)=\{u_2\}$.  Hence $uv\in E_r$, otherwise, we obtain a bad $2$-coloring of $G+wu$ from $c$ by coloring $wu$ red and then recoloring all edges incident with $w$ in $D_2$ blue. Therefore $v\in R_2$.\proofsquare

By symmetry, for any     $w \in R_1$,  if $w$ is adjacent to exactly one vertex, say $v$, in   $H\less P$, then $v\in Q_2 $. We next count the number of edges in $H$. 
Since $N_r(u_1)\cap N_r(u_2)=\emptyset$, it follows from  \pr{e:rednbr} that  for each $w\in P$, $e_H(w, Q\cup R) \ge2$ and so $e_H(P, Q\cup R)\ge 2|P|$. 
Let $Q_1^*$ be the set of vertices  $w\in Q_1$ such that $w$ is adjacent to exactly one vertex in  $H\less P$. Similarly, let $R_1^*$ be the set of vertices  $w\in R_1$ such that $w$ is adjacent to exactly one vertex  in  $H\less P$.   By \pr{e:1rednbr},  $e_H(Q_1^*, R_2)\ge |Q_1^*|$ and $e_H(R_1^*, Q_2)\ge |R_1^*|$. Notice that for any $w\in (Q_1\cup R_1)\less (Q_1^*\cup R_1^*)$, $w$ is adjacent to at least two vertices in $H\less (P\cup Q_1^*\cup R_1^*)$ and so 
$e(H\less (P\cup Q_1^*\cup R_1^*))\ge |Q_1\less Q_1^*|+|R_1\less R_1^*|=|Q_1| + |R_1|-|Q_1^*|-|R_1^*|$. Therefore 
\begin{align*}
e(H)&=e_H(P, Q\cup R)+e_H(Q_1^*, R_2)+e_H(R_1^*, Q_2)+e(H\less (P\cup Q_1^*\cup R_1^*))\\
&\ge 2|P|+|Q_1^*|+|R_1^*|+|Q_1| + |R_1|-|Q_1^*|-|R_1^*|\\
&=2|P|+|Q_1|+|R_1|.
\end{align*}
 Note that $e_G(V(D_1), Q\cup R)\ge |Q_1|+2|Q_2|+|R_1|+2|R_2|=|Q|+|R|+|Q_2|+|R_2|$. We see that $e(H)+e_G(V(D_1), Q\cup R) \ge (2|P|+|Q_1|+|R_1|)+(|Q|+|R|+|Q_2|+|R_2|)=2(|P|+|Q|+|R|)\ge 2n-2k+2$. 
  By \pr{e:edgecount}, 
\begin{align*}
e(G) &\ge e(H)+e_G(V(D_1), Q\cup R)  + \sum_{i = 1}^p e(G[V(D_i)]) \\
&\ge (2n-2k+2)+ \left(\frac{1}{2} \left\lceil \frac{k}{2} \right\rceil - \frac{1}{2} \right) n -\left(\frac{1}{2} \left\lceil \frac{k}{2} \right\rceil - \frac{1}{2} \right) k\\
 &= \left( \frac{3}{2} + \frac{1}{2} \left\lceil \frac{k}{2} \right\rceil \right) n - \left(\frac{1}{2} \left\lceil \frac{k}{2} \right\rceil + \frac{3}{2} \right) k+2\\
 &= \left( \frac{3}{2} + \frac{1}{2} \left\lceil \frac{k}{2} \right\rceil \right) n -c
\end{align*}
where $c=\left(\frac{1}{2} \left\lceil \frac{k}{2} \right\rceil + \frac{3}{2} \right) k -2$. 
\medskip

This completes the proof of  Theorem~\ref{K3Tk}.
\proofsquare

\noindent {\bf Conclusion}.   For the graphs $G_{odd}$ and $G_{even}$ in the proof of Theorem~\ref{K3T4}, we want to point out here that we found the graph  $G_{odd}$ when $d_b(y)=1$, $d_b(z)=2$,  and $G_r=J$ with $|B|=2$ and $|C|=4$; and the graph  $G_{even}$ when $d_b(y)=d_b(z)=2$,  and  $G_r=J$ with $|B|=3$ and $|C|=2$. We believe that  the method we developed in this paper can be applied to determine $sat(n, \mathcal{R}_{\min}(K_p, T_k))$ for any given tree $T_k$ and any $p\ge3$.

\section*{Acknowledgments}

The authors would like to thank Christian Bosse, Michael Ferrara, and Jingmei Zhang for their helpful discussion.  The authors thank the  referees for helpful comments.
\medskip

\end{document}